\documentclass[english,onecolumn]{IEEEtran}
\usepackage{graphicx} 
\usepackage{geometry}
\usepackage{amsthm}
\usepackage[T1]{fontenc}
\usepackage{amsmath}
\usepackage{graphicx}
\usepackage{amssymb,bbm}
\usepackage{mathtools}
\usepackage{float}

\usepackage{xcolor}

\usepackage{mathtools}

\newtheorem{lemma}{Lemma}
\newtheorem{problem}{Problem}
\newtheorem{example}{Example}
\newtheorem{proposition}{Proposition}
\newtheorem{theorem}{Theorem}
\newtheorem{definition}{Definition}
\newtheorem{assumption}{Assumption}
\newtheorem{remark}{Remark}

\newcommand{\zeros}[1]{0_{#1}}
\newcommand{\eye}[1]{I_{#1}}
\newcommand{\STM}[2]{\Phi(#2,#1)}

\title{A Hybrid Perspective on Suboptimal Mixed-Integer Quadratic Programming}

\author{Luke Fina$^{*}$, Christopher Petersen$^{*}$\thanks{$^{*}$Department of Mechanical and Aerospace Engineering, University
of Florida, Gainesville, FL USA. Emails: \textit{l.fina@ufl.edu,~c.petersen1@ufl.edu. }}}
\begin{document}

\maketitle
\begin{abstract}
    This work solves suboptimal mixed-integer quadratic programs recursively for feedback control of dynamical systems. The proposed framework leverages parametric mixed-integer quadratic programming (MIQP) and hybrid systems theory to model a recursive MIQP feedback controller and a feedback controlled dynamical system. The proposed hybrid framework theoretically encodes the suboptimal part via solver parameters as bounded perturbations from the optimal solution set. The stability of the proposed hybrid framework is theoretically guaranteed and validated through MIQP model predictive control simulations with multiple solver parameters.
\end{abstract}

\maketitle

\section{Introduction} 
Optimization-based feedback control is dependent on recursive optimization algorithms updating as new information from a dynamical system arrives, like position and disturbances updates. Approximation error is well known to occur for feedback control and classically is viewed as uncertainty in the plant model \cite{doyle2013feedback}. In recent years, researchers have begun to analyze approximation error bounds in the controller due to optimization algorithm performance in optimization-based feedback control \cite{hauswirth2024optimization}. This work extends this approximation error bound analysis to the field of recursive mixed-integer optimization and analyzes the asymptotic behavior of the mixed-integer feedback control of a dynamical system, where the term mixed-integer means there are discrete and continuous decision variables.

An optimization-based feedback controller with algorithm-based approximation error analysis can be viewed as a discrete-time dynamical system \cite{hauswirth2024optimization}. For example, treating a continuous optimization algorithm as a discrete-time dynamical system recently provided new insights into algorithm phenomena that have previously eluded theoretical understanding for the accelerated Nestrov method \cite{muehlebach2019dynamical}. Treating mixed-integer algorithms as discrete-time dynamical systems is less explored \cite{hauswirth2024optimization}. One work that explores this perspective, solves integer optimization as a collection of coupled oscillators \cite{guanchun2023dynamical}. Aside from this work, the discrete-time dynamical systems perspective for mixed-integer algorithms has seldom been explored to the author's knowledge.

One important distinction to make is that  recursive continuous optimization algorithms have analytical convergence rates that allow for strong theoretical discrete-time dynamics analysis \cite{bertsekas1997nonlinear}, while mixed-integer optimization algorithms do not have theoretically derivable convergence rates \cite{pia2017mixed}. However, modern mixed-integer solvers combine a collection of heuristics, branch-and-bound, and cuts to achieve fast solutions in practice \cite{clautiaux2024last,pia2017mixed}. 
To that end, this work abstracts mixed-integer algorithms to a black box solver, introduces relevant black box solver terms, and analyzes the asymptotic behavior of the mixed-integer solver and feedback controlled dynamical system.

Quadratic programs (QPs) are prevalent in recursive continuous optimization for feedback control, e.g. model predictive control and control barrier functions \cite{ames2019cbf, morari1999model}. Due to the prevalence of quadratic programming in feedback control, our work focuses on the mixed-integer version of quadratic programming, mixed-integer quadratic programming (MIQP). MIQPs can model control scenarios that QPs cannot, e.g. when the control is integer in model predictive control \cite{mcallister2022advances} or signal temporal logic \cite{kurtz2022mixed}. This work focuses on recursive MIQPs, but the theory can be easily extended to recursive mixed-integer linear programming and recursive indefinite quadratic programming. 

To be precise, this work leverages the theory of parametric MIQP for recursive MIQPs to prove the MIQP solver feedback controller and feedback controlled dynamical system with sufficient regularity conditions is a hybrid system. A blend of continuous time and discrete time systems can be theoretically described by hybrid systems \cite{goebel2009hybrid,sanfelice2021hybrid}. The MIQP solver operating on hardware is modeled over discrete time and the feedback controlled dynamical system is modeled over continuous time. Recently, continuous optimization algorithm works have explored a hybrid system framework \cite{hendrickson2025distributed, hustig2024uniting}. Finally, this work explores specific MIQP solver parameters useful for practitioners.

To summarize, the main contributions of this paper are the following:
\begin{itemize}
    \item A novel hybrid systems framework for recursive MIQPs is proposed, and the framework generalizes to recursive mixed-integer linear programs and recursive indefinite quadratic programs by leveraging parametric MIQPs.
    \item Sufficient regularity conditions are presented for the recursive MIQP and feedback controlled dynamical system to satisfy the basic hybrid conditions (Theorem~\ref{thm:hybrid_cond_MIQP}).
    \item Stability conditions are derived for solutions to the proposed hybrid system with suboptimal solver parameters  (Theorem~\ref{thm:perturb}).
    \item Useful suboptimal solver parameters are explored by solving MIQP model predictive control in simulation.
\end{itemize}

The rest of the paper is organized as follows. Section~\ref{sec:prelims} presents preliminaries for the recursive MIQP and the hybrid system formulation. Section~\ref{sec:solver_processes} defines solver processes and solver parameters. Section~\ref{sec:stability}  presents stability guarantees for the proposed hybrid framework. Section~\ref{sec:simulations} explores suboptimal solver parameters for MIQP model predictive control and the impact on hybrid solutions. Section~\ref{sec:conclusion} concludes. Section~\ref{sec:appendix} includes technical details necessary for the proofs in this work.

\textbf{Notation and Terminology:}
In this work, parametric optimization is leveraged to theoretically analyze recursive optimization. Parametric optimization is based on point set mappings. Point to set mappings are used within the optimization literature and are denoted as the power set $2^{\mathbb{R}^{p}},$ which is the family of all sets \cite{kharazishvili1998set}. Set value mappings by definition are mappings from an arbitrary set $X$ to the power set.  This work treats point to set mappings as the more general set value mappings throughout this work \cite{rockafellar2009variational}[Ch 5].
A set valued map $M$ from $S\subseteq \mathbb{R}^{p}$ to the power set of some Euclidean space $\mathbb{R}^{p}$ is denoted by $M:S\rightrightarrows\mathbb{R}^{p}.$ The domain $\text{dom}M$ of a set value map is defined as $\text{dom}M=\{x\in \mathbb{R}^{p} \vert M(x) \neq \emptyset\}.$  Set value mappings are also the fundamental mathematical tool for differential inclusions. Discrete time variants of differential inclusions are termed difference inclusions and are used for explicitly modeling the optimization problem in the proposed hybrid systems framework. A function $\alpha: \mathbb{R}_{\geq 0} \to \mathbb{R}_{\geq 0}$ is a class $\mathcal{K}_{\infty}$ function if $\alpha$ is zero at zero, continuous, strictly increasing, and unbounded. A function $\beta:\mathbb{R}_{\geq 0}\times \mathbb{R}_{\geq 0} \to \mathbb{R}_{\geq 0}$ is a class $\mathcal{KL}$ function if it is non-decreasing in its first argument, non-increasing in its second arguments, i.e, $\lim_{r\to 0^{+}}\beta(r,s)=0$ for each $s\in \mathbb{R}_{\geq 0},$ and $\lim_{s\to \infty}\beta(r,s)=0$ for each $r\in \mathbb{R}_{\geq 0}.$ The notation $\vert \cdot \vert_{\mathcal{A}}$ denotes the Euclidean norm for a point $\cdot$ to the set $\mathcal{A}.$ This work denotes the real numbers as $\mathbb{R},$ the rational numbers as $\mathbb{Q},$ and $x+\mathbb{B}$ denotes a closed Euclidean unit ball center at $x.$ 

\section{Preliminaries and Problem Statement}\label{sec:prelims} 
This section connects notions of regularity from a parametric MIQP and hybrid systems context, then formulates a hybrid systems framework for parametric MIQPs with sampled continuous dynamics changing over time. Section~\ref{sec:stability} derives stability conditions for the proposed hybrid system.

\subsection{Parametric MIQPs Preliminaries}
 MIQPs in feedback control are solved recursively, we leverage the theoretical framework of parametric optimization to theoretically analyze the recursive MIQP. Varying information can be constructed in a constraint or objective term at every sample time $k$. In this work, the parametric MIQP focuses on information varying in a linear objective term and affine constraints, where the parameter is a sampled state from a continuous dynamical system. The connection between parametric optimization and time-varying convex optimization literature inspires the connection between recursive MIQPs and parametric optimization in this work \cite{simonetto2020time}. 
 This work studies a parametric MIQP of the following form
 
\begin{problem}\label{prob:MIQP}
 \begin{align*}
    V(c(k),\eta(k))= &\underset{y}{\text{minimize}}  \enskip y^{T}Qy+c(k)^{T}y\\
    &\text{subject to} 
    \\
    &\qquad Ay \leq \eta(k)
    \end{align*}
where  $y_{1},...,y_{s}\in \mathbb{Z}^{s}$ are the discrete variables that are part of the full vector $y \in \mathbb{R}^{n}$, $c(k)\in\mathbb{R}^{n},~\eta(k)\in\mathbb{R}^{m}$ are arbitrary parameter vectors, $A\in \mathbb{R}^{m\times n},~Q\in \mathbb{R}^{n\times n}.$ 
\end{problem}
The terms $c(k)$ and $\eta(k)$ vary with new information at constant sample times $k$, which will take an explicit form in the hybrid systems formulation. Information invariant terms $Ay\leq b$ and $d^{T}y$ can also be included in the formulation of Problem~\ref{prob:MIQP}, but for notational simplicity these terms are omitted.

We start with regularity for the solutions to Problem~\ref{prob:MIQP}, because regularity provides stability conditions for the solutions, i.e. conditions where small changes in the parameters $c(k)$ and $\eta(k)$ mildly perturb solutions to Problem~\ref{prob:MIQP}. For MIQPs this can be accomplished with rational matrices. 
\begin{assumption}\label{assum:rational}
The objective function matrix $Q$ and constraint matrix $A$ are rational matrices, i.e.$~Q\in \mathbb{Q}^{n\times n}$ and $~A\in \mathbb{Q}^{m\times n}$.
\end{assumption}
The above assumption is necessary to prove stability of solving the parametric mixed-integer quadratic program \cite{bank1984stability}, which is reasonable in practice \cite{gautschi2011numerical}.

Problem~\ref{prob:MIQP} can be described more generally as multiple set value mappings. For notational simplicity, this work denotes the objective function as $\phi(y,c(k)):= y^{T}Qy+c(k)^{T}y.$ Let $H(\eta(k)): \mathbb{R}^{m}\rightrightarrows\mathbb{R}^{n},$ denote a feasible set value mapping of the constraints. Let   $V(c(k),\eta(k)): \mathbb{R}^{n}\times \mathbb{R}^{m}\rightrightarrows\mathbb{R}$ denote a set value mapping of the value function, where the value function is the objective function evaluated at the minimizer $y^{*},$ i.e $\phi(y^{*},c(k)).$ Then construct 
\begin{equation}\label{eq:suboptimal_setvalue_mapping}
\psi_{\epsilon}(c(k),\eta(k),\epsilon):=\{ y\in \mathbb{R}^{n} \vert \enskip y \in H(\eta(k)), \phi_{\epsilon}(y,c(k))\leq V(c(k),\eta(k))+\epsilon \}
\end{equation}
as the $\epsilon$-optimal set value mapping, $\psi_{\epsilon}: \mathbb{R}^{n} \times \mathbb{R}^{m} \rightrightarrows\mathbb{R}^{n}$ \cite{bank1982non,bank1984stability}, where $\epsilon\in\mathbb{R}$ denotes suboptimality.  Next, lemmas are presented that describe the regularity of $\psi_{\epsilon}(c(k),\eta(k),\epsilon)$ with varying $c(k)$ and $\eta(k)$ when Assumption~\ref{assum:rational} holds.

\begin{lemma}\label{lem:subopt_mapping_usc}Let Assumption~\ref{assum:rational} hold, then the $\epsilon$-optimal set value mapping defined in Eq (\ref{eq:suboptimal_setvalue_mapping}) is upper semi-continuous. 
\end{lemma}
See Appendix~\ref{appen:MIQP_stability} for technical details on proving Lemma~\ref{lem:subopt_mapping_usc}.
The upper semi-continuity property holding for Eq.~(\ref{eq:suboptimal_setvalue_mapping}) is necessary to prove Problem~\ref{prob:MIQP} satisfies basic conditions for a well behaved hybrid system, because the decision variable $y$ will be a state in the constructed hybrid system.

\subsection{Hybrid System Preliminaries}
Regularity for hybrid systems will be defined next before introducing the proposed hybrid system. Regularity for a hybrid system will provide conditions for well posedness of the hybrid solutions \cite{goebel2009hybrid}, which in this work is relevant to describe the MIQP-based feedback controller and the feedback controlled dynamical system. Well-posedness is necessary for proving robustness of the hybrid system in Theorem~\ref{thm:perturb}. 

We denote a hybrid solution as $\pi(t,j)$, where $t$ are flow times (i.e., time arcs where the system flows continuously) and $j$ are jump times (i.e., discontinuities where the hybrid system jumps states) with a hybrid state $\chi(t,j)$ denoted by $\chi$ for notational brevity. A general hybrid system is defined as $\mathcal{H}=(C,F,D,G)$ with

\begin{align*}
&\chi \in C, \quad \dot{\chi}= F(\chi),\\
&\chi \in D, \quad \chi^{+}=G(\chi),
\end{align*}
where $\chi\in \mathbb{R}^{h}$ denotes the hybrid system state, $F$ denotes the set value flow map, $C$ denotes the flow set, $G$ denotes the set value jump map, $D$ denotes the jump set \cite{goebel2009hybrid}[Theorem 6.8], and $\dot{\chi},~\chi^{+}$ denote the evolution of the hybrid state in the flow and jump map respectively. This work models the MIQP solver and dynamical system in discrete and continuous time respectively. In hybrid system terminology this means the set value mapping for Problem~\ref{prob:MIQP}, i.e $y\in\psi_{\epsilon}$, jumps and the continuous dynamical system $\dot{x}=f(x,u)$  flows. Regularity for hybrid systems can be defined by the hybrid basic conditions,
\begin{definition} (Hybrid Basic Conditions)\label{def:hybrid_conds}
\cite{goebel2009hybrid}
\begin{enumerate}
    \item $C$ and $D$ are closed subsets of $\mathbb{R}^{h}$;
    \item The flow map $F$ is outer semi-continuous and locally bounded relative to $C$, $C\subset \text{dom} F$, and $F(\chi)$ is convex for every $\chi\in C$;
    \item The jump map $G$ is outer semi-continuous and locally bounded relative to $D$, $D\subset \text{dom} G$.
\end{enumerate}
\end{definition}

The set value mapping $y\in\psi_{\epsilon}(c(k),\eta(k),\epsilon)$ from Problem~\ref{prob:MIQP} combine with continuous dynamics $\dot{x}=f(x,u)$ in a sample and hold paradigm to form the proposed hybrid system, $\mathcal{H}_{\text{MIQP}}.$ Sampling $\dot{x}=f(x,u)$ gives an explicit form for the generic parameter vectors $\eta(k)$ and $c(k)$, which are $\eta(k):=x(\tau)$ and $c(k):=\eta(x(\tau)).$ The term $x(\tau)$ denotes sampled continuous dynamics. The term $\eta(x(\tau))$ denotes a time-varying objective term that evolves based on the sampled dynamics. At a high level the $\epsilon$ term encodes suboptimality from solver parameters by treating suboptimality as a perturbation from the optimal solution set and is defined further in Section~\ref{sec:solver_processes}. 

For brevity, if a state is not explicitly defined in the flow or jump map, it has no change. The proposed hybrid system $\mathcal{H}_{\text{MIQP}}$ jump map is

\begin{equation}\label{eq:G_MIQP}
    G(\chi)= \begin{bmatrix}y^{+} \in \psi_{\epsilon}\big(x(\tau),\eta(x(\tau)),\epsilon \big) \\

    \tau^{+}=0
   \end{bmatrix},
    \end{equation}
    where $\psi_{\epsilon}$ is defined in Eq.~(\ref{eq:suboptimal_setvalue_mapping})
     and the jump set is
    \begin{equation}\label{eq:D_MIQP}
    D= \{y \in \mathbb{R}^{n}, ~\tau=\Delta t\},
    \end{equation}
    where $\tau$ is a timer and $\Delta t$ denotes sampling time.
   The flow map is
    \begin{equation}\label{eq:F_MIQP}
    F(\chi) = \begin{bmatrix} \dot{x}=f(x,\kappa(y)) \\ \dot{\tau}=1\end{bmatrix},
    \end{equation}
    where $\kappa(y)$ is a function to extract the control portion of the decision vector $y$ from Problem~\ref{prob:MIQP}. The flow set is
    \begin{equation}\label{eq:C_MIQP}
    C= \{\tau \in \mathbb{R}, \enskip\tau \in  [0,\Delta t],~x\in \mathbb{R}^{d}\},
    \end{equation}
    and the hybrid state is
     \begin{equation}\label{eq:hybrid_states}\chi=(y,x,\tau).
     \end{equation}

 The timer logic is as follows: increment towards $\Delta t$ until $\tau=\Delta t$ during flows, then the hybrid system jumps and the timer is held constant, the timer resets to zero, and repeats. 

  \subsection{Connecting Parametric MIQPs and Hybrid Systems}
  This section has introduced two notions of regularity: 1) Assumption~\ref{assum:rational} and 2) Definition~\ref{def:hybrid_conds}. Our first result proves, under sufficient conditions, the two notions are equivalent allowing us to leverage hybrid systems theory for the remainder of this work. 
  
  Before the result, another preliminary lemma connects outer semi-continuity from Definition~\ref{def:hybrid_conds} and upper semicontinuity from $\psi_{\epsilon}\big(x(\tau),\eta(x(\tau)),\epsilon\big).$

\begin{lemma}\cite{goebel2009hybrid}
\label{lem:osc_usc}
(Outer semi-continuous vs upper semi-continuity) Let $M$ be a set-valued mapping. Consider $\chi$ such that $M(\chi)$ is closed. If $M$ is upper semi-continuous at $\chi$, that is, every $\vartheta>0$ there exists $\delta>0$ such that $\chi'\in \chi+\delta \mathbb{B}$ implies $M(\chi')\subset M(\chi)+\vartheta\mathbb{B},$ then $M$ is outer semi-continuous at $\chi$. If $M$ is locally bounded at $\chi$, then the reverse implication is true.
\end{lemma}

For brevity, this work will assume the dynamics satisfy the basic hybrid conditions, which are satisfied by many nonlinear continuous systems \cite{goebel2009hybrid}.
\begin{assumption}\label{assum:hybrid_conds_cont_dyn}The dynamics $\dot{x}=f(x,\kappa(y))$ are outer semi-continuous, locally bounded relative to $C,$  $f(x,\kappa(y))$ is convex for every $x \in C.$
\end{assumption}

Using Lemma~\ref{lem:osc_usc} under Assumption~\ref{assum:rational} and Assumption~\ref{assum:hybrid_conds_cont_dyn}, the next result guarantees Definition~\ref{def:hybrid_conds} is satisfied for optimal solutions of $\mathcal{H}_{\text{MIQP}}$.
\begin{theorem}\label{thm:hybrid_cond_MIQP}Let Assumption~\ref{assum:rational} and Assumption~\ref{assum:hybrid_conds_cont_dyn}  hold. Suppose Problem~\ref{prob:MIQP} is optimal, i.e. $\epsilon=0.$ Then $\mathcal{H}_{\text{MIQP}}$ satisfies the hybrid basic conditions.
\end{theorem}
\begin{IEEEproof}
By construction, the sets $C$ , Eq. (\ref{eq:C_MIQP}), and $D$ , Eq. (\ref{eq:D_MIQP}) are subsets of $\mathbb{R}^{h}.$ Therefore condition 1 of Definition~\ref{def:hybrid_conds} that $C$ and $D$ are closed subsets of $\mathbb{R}^{h}$ is true.

Condition 2 of Definition~\ref{def:hybrid_conds} is that (i) the flow map is outer semi-continuous and locally bounded relative to the flow set, (ii) the flow set is a subset of the domain of the flow map, and (iii) the flow map is convex for every state in the flow set.

The timer is a constant function, which is a continuous function. Therefore, the timer satisfies condition 2 of the Hybrid Basic Conditions. By Assumption~\ref{assum:hybrid_conds_cont_dyn}, the dynamics $\dot{x}=f(x,\kappa(y))$ are outer semi-continuous and locally bounded relative to $C$ and $f(x,\kappa(y))$ is convex for every $x\in C,$ then condition 2 of Definition~\ref{def:hybrid_conds} is satisfied.

Condition 3 of Definition~\ref{def:hybrid_conds} is that the jump map, Eq. (\ref{eq:G_MIQP}), needs to be proven to be (a) outer semi-continuous, (b) locally bounded relative to the jump set, Eq. (\ref{eq:D_MIQP}), and (c) the jump set is a subset of the domain of the jump map.

Using Lemma~\ref{lem:subopt_mapping_usc}, the term $\psi_{\epsilon}(x(\tau),\eta(x(\tau)),0)$ is upper semi-continuous. Condition (a) holds using Lemma~\ref{lem:osc_usc}, because upper semi-continuity implies outer semi-continuity. 


Condition (b) holds because upper semi-continuity by  definition is locally bounded (See Definition~\ref{def:u_s_c}), i.e $\psi_{\epsilon}(x(\tau),\eta(x(\tau)),0)$ is locally bounded for the domain $y\in \mathbb{R}^{n},$ and the timer $\tau$ is trivially within the domain because it is held constant in the jump set, Eq. (\ref{eq:D_MIQP}).


To prove $D\subset \text{dom}G,$ by definition 
$
\text{dom}G = \{ \chi \in \mathbb{R}^{h} \vert G(\chi) \neq \emptyset\},
$
which can also equivalently be rewritten as a new set value mapping
\[
G_{D}(\chi) = \begin{cases}
    G(\chi) \enskip \text{if } \chi\in D\\
    \emptyset \enskip \text{if }\chi\notin D.
\end{cases}
\]
Therefore,
$G_{D} = \text{dom}G$ and by construction
$
D\in G_{D}\implies D \subset \text{dom}G,
$ i.e. Condition (c) holds. Then all three conditions for the Basic Hybrid Conditions in Definition~\ref{def:hybrid_conds} hold. This completes the proof.
\end{IEEEproof}
The above result allows this work to treat the feedback control system of Problem~\ref{prob:MIQP} and $\dot{x}=f(x,\kappa(y))$  as a typical hybrid system, $\mathcal{H}_{\text{MIQP}}$. We use hybrid systems tools later to derive stability conditions for solutions to $\mathcal{H}_{\text{MIQP}}$.

\subsection{Generality of Approach}

Quadratic programming is employed often for recursive optimization in feedback control \cite{ames2019cbf,morari1999model}. This motivates focusing on the mixed-integer counter part MIQPs, however the results in this work are not limited to only MIQPs.  The hybrid system $\mathcal{H}_{\text{MIQP}}$ is quite general and can model alternative time-varying optimization problems by leveraging parametric optimization theory, e.g. parametric indefinite quadratic programs and parametric mixed-integer linear programs. This generality holds due to only assuming Assumption \ref{assum:rational} for Problem \ref{prob:MIQP} without any compactness assumptions.
\begin{remark}
    These results also apply for indefinite continuous quadratic programs without needing Assumption~\ref{assum:rational} and for mixed-integer linear programs without needing the quadratic objective function term. 
\end{remark}

In addition to general theoretics of Problem~\ref{prob:MIQP}, the next section focuses on generic MIQP solver properties to solve Problem~\ref{prob:MIQP} before diving into stability of solutions to the hybrid system $\mathcal{H}_{\text{MIQP}}$ in Section~\ref{sec:stability}

\section{Solver Processes}\label{sec:solver_processes} This section introduces definitions for generic MIQP solver progress and solver parameters that are used further in Section~\ref{sec:simulations}. Revisiting the term $\epsilon$ in Eq. (\ref{eq:suboptimal_setvalue_mapping}), for suboptimal continuous optimization, $\epsilon$ is based on a theoretical convergence rate and step size of particular algorithms \cite{amiri2025practical}. In contrast, this work generalizes these concepts for black box MIQP solvers when theoretical convergence rates are not available.

\subsection{Solver Iterates}
In continuous optimization the notion of algorithm improvement towards the optimal solution is described by iterations. Iterations are ill defined for mixed-integer solvers, because
commercial and open-source mixed-integer solvers, like Gurobi, SCIP, and HiGHS, use a combination of heuristics, cutting planes, and branching \cite{conforti2014integer,gurobi_MIP,SCIP_MIP,highs_MIP}. Therefore, this section defines a solver iterate.

Let us view any optimization solver as a process with inner and outer loops. For a gradient descent algorithm, the inner and outer loop are equivalent. A gradient descent algorithm with an adaptive step size selection process has an inner and outer loop. The inner loop is step size selection and the outer loop is step size selection plus a step of the gradient. Mixed-integer programming solvers are significantly more complex than adaptive gradient descent. Our next definition provides a notion of solver progress applicable for mixed-integer solvers.
\begin{definition}\label{def:solver_compute}A solver iterate is any outer loop within the solver that starts with an initial guess $y_{0}$ and terminates after $i_{f}$ outer loops under an exit criteria.
 \end{definition}
 Suppose, again, a gradient descent algorithm, then a solver iterate is equivalent to an iteration of gradient descent, since inner and outer loops are equivalent. For gradient descent with an adaptive step size selection, picking a new step size and one gradient descent step would qualify as a solver iterate.  
 Next, let us illustrate a case for mixed-integer solvers. In the mixed-integer linear program (MILP) solvers HiGHs and CPLEX, one inner loop by default is a warm start method termed crash start \cite{galabova2023presolve}. Crash start finds a partially feasible solution to initialize the main MILP solver. One inner loop of crash start and loop of the main solver would qualify as a solver iterate. We make this distinct because solving a mixed-integer problem typically involve solving continuous relaxations within branching-and-bounding the integer variables. This leads to iterates for continuous algorithms that are inner loops and solver iterates for the overall solver, which are outer loops. 
 

\subsection{Solver parameters}\label{sec:solver_params}
Next, solver parameters are split into two classes: existence solver parameters and suboptimal solver parameters.

\begin{definition}
    Existence solver parameters impact if a numerical solution exists.
\end{definition}
Examples of existence solver parameters are hard time limits, hard memory limits, and numeric stability. MIQP solvers run on hardware and that hardware has its own limitations \cite{gurobi_MIP}. Therefore, existence solver parameters may vary computationally between hardware devices. Existence solver parameters interrupt the solver iterate process and exit with an error. To illustrate, if the hard memory limit on the hardware is exceeded, then the solver exits with an error. From a hybrid systems context, existence solver parameters can lead to a maximal solution of the hybrid system. This class of solver parameters will be explored in a future work to explore the experimental relationship between compute hardware and existence solver parameters.

In contrast, solver parameters that impact suboptimality of the final solution are introduced next and are the focus of the rest of this work.

\begin{definition}\label{def:subopt_hyper_param}
    Suboptimal solver parameters impact how close the final solution at $i_{f}$ solver iterates is to the optimal solution.
\end{definition}

Three useful exit criteria based suboptimal solver parameters are the following: 1) soft time limit, 2) soft memory limit, and 3) iteration limit. Soft time limit sets an exit criteria to limit the MIQP solver run time. Soft memory limit  sets an exit criteria to limit the memory allocated to the solver to find a solution. The term soft denotes the solver parameter exits with no error in contrast to the hard solver parameters. Conceptually the solver exits close to the exact exit criteria, but will wait to prevent outputting an error. While solving a MIQP, quadratic programming relaxations are solved in a branching and bounding fashion. Iteration limit  sets an exit criteria to limit the number of quadratic programming iterations. In Gurobi the iteration limit parameter limits simplex iterations, which is likely due to internal reformulations of MIQPs to mixed-integer linear formulations. More details can be found in \cite{gurobi_MIP}.

 Suboptimal solver parameters in Eq. (\ref{eq:suboptimal_setvalue_mapping}) are modeled as $\epsilon$, and in a hybrid systems context perturb solutions from the optimal solution set which is discussed in Section~\ref{sec:stability}. 
 The three exit criteria solver parameters are suboptimal as defined in Definition~\ref{def:subopt_hyper_param} and are validated in simulation in Section~\ref{sec:simulations}.

\section{Stability Analysis}\label{sec:stability}
This section addresses when the solutions to the hybrid system with potentially suboptimal solver parameters are stable. We employ warm starting Problem~\ref{prob:MIQP} every sample time to ensure stability and robustness without an analytical solver convergence rate for the hybrid system, $\mathcal{H}_{\text{MIQP}}$. Warm starting is a general term that entails using the solution from the last solver iterate at sample time $k$ as the initial guess in the next solver iterates at sample time $k+1$.

Let us define a stabilizing set for the hybrid system $\mathcal{H}_{\text{MIQP}}$ as
\[
\mathcal{A} =  Y\times \{ 0\} \times [0,\Delta t],
\]
where $Y$ is a closed set of solutions for Problem~\ref{prob:MIQP} defined by the $\epsilon$-optimal set value mapping $\psi_{\epsilon}$ where $\epsilon=0$ for all sample times,  $\{0\}$ is a compact set of the equilibrium for the dynamics $\dot{x}=f(x,\kappa(y))$, and $[0,\Delta t]$ is the compact set of timer values.  Conceptually this is the set of all optimal solutions that stabilize the dynamics within the allotted sample time.

A Lyapunov function can be defined for the above set as
 \begin{equation}\label{eq:lyap_fxn}
 L = (\phi_{\epsilon}-V)+\frac{1}{2}x^{2},
 \end{equation}
 where the term $(\phi_{\epsilon}-V)$ encodes deviations from optimality of Problem~\ref{prob:MIQP} and the term $\frac{1}{2}$ encodes deviation from the equilibrium.
The upcoming asymptotic stability analysis leverages feasible warm starting with an optimization solver (e.g. when solvers are given a ``good" initial guess) and the Lyapunov function extensively to prove stability of the solutions to the  $\mathcal{H}_{\text{MIQP}}.$ Prior to the next result, an assumption on how warm starting relates to the hybrid system $\mathcal{H}_{\text{MIQP}}$ is established.
\begin{assumption}\label{assum:warm_start}
Warm starting with an initial feasible guess $\mathcal{H}_{\text{MIQP}}$ implies
 there exists $K_{\infty}$ functions $\alpha_{1},~\alpha_{2}$ such that  $\alpha_{1}(\vert \chi \vert_{\mathcal{A}})\leq L(\chi) \leq \alpha_{2}(\vert \chi \vert_{\mathcal{A}})$ holds.

\end{assumption}
The above assumption is a technical condition meaning the Lyapunov function in Eq. (\ref{eq:lyap_fxn}) is bounded by smooth functions, which is necessary for the next major result.

Before presenting the next major result, an assumption on the relationship between the control $\kappa(y)$ from Problem~\ref{prob:MIQP} and the dynamical system  $\dot{x}=f(x,\kappa(y))$ is needed. 
\begin{assumption}\label{assum:stable_control} The control law $\kappa(y)$ for $\dot{x}=f(x,\kappa(y))$ is asymptotically stable.
\end{assumption}
Feasible optimization-based control MIQPs satisfy this condition, e.g. warm started suboptimal model predictive control \cite{scokaert1999suboptimal}. 
We next construct sufficient conditions under which solutions to the set $\mathcal{A}$ is uniform globally pre-asymptotically stable for $\mathcal{A}$ for $\mathcal{H}_{\text{MIQP}}.$ 

\begin{lemma}\label{lem:asym_stable} Let Assumption~\ref{assum:warm_start} and Assumption~\ref{assum:stable_control} hold. Suppose the hybrid system $\mathcal{H}_{\text{MIQP}}$ is warm started with an initial feasible guess, then any solution $\pi(t,j)$ to the set $\mathcal{A}$ is uniformly globally pre-asymptotically stable when solving Problem~\ref{prob:MIQP} to optimality.
\end{lemma}
\begin{IEEEproof} 
By \cite{sanfelice2021hybrid}[Theorem 3.19 3.a], the set $\mathcal{A}$ is uniformly global pre-asymptotic stability for $\mathcal{H}_{\text{MIQP}}$ if $\mathcal{A}$ is closed and the following holds:
 there exists class $K_{\infty}$ functions $\alpha_{1},\alpha_{2}$ such that $\alpha_{1}(\vert \chi \vert_{\mathcal{A}})\leq L(\chi) \leq \alpha_{2}(\vert \chi \vert_{\mathcal{A}})$ holds, which is true under Assumption~\ref{assum:warm_start};
     the flow dynamics are asymptotically stable, i.e. $\dot{L}(\chi)\leq -\rho (\vert \chi\vert_{\mathcal{A}})~\forall\chi\in C,$  which is true under Assumption~\ref{assum:stable_control};
     the jump dynamics are asymptotically stable, i.e. $\Delta L(\chi) \leq -\rho(\vert \chi \vert _{\mathcal{A}}),$ which is true when Problem~\ref{prob:MIQP} is solved to the optimal solution meaning $\phi_{\epsilon}-V=0$.
This completes the proof.
\end{IEEEproof}

\subsection{Suboptimal Solver Parameters Stability Analysis}
 This subsection derives stability for the influence of suboptimal solver parameters on solutions to $\mathcal{H}_{\text{MIQP}}$, where the influence perturbs solutions away from the set $\mathcal{A}.$  
 
 The suboptimal solver parameters perturbing on the hybrid system $\mathcal{H}_{\text{MIQP}}$ form a new hybrid system $\mathcal{H}_{e}$ defined as

\begin{equation*}
    G_{e}(\chi):= \begin{bmatrix}y^{+} \in \psi_{\epsilon}\big(x(\tau),\eta(x(\tau)),e_{y} \big) \\

    \tau^{+}=0
   \end{bmatrix},
    \end{equation*}
\begin{equation*}
    F_{e}(\chi) = \begin{bmatrix} \dot{x}=f(x,\kappa(y)) \\ \dot{\tau}=1+e_{\tau}\end{bmatrix},
    \end{equation*}
where  $e_{\tau}=[-1,\infty)$ models the time perturbation from the solver parameters, and $e_{y}=(0,\hat{e}_{y}]$ models the suboptimality from the solver parameters, where $e_{y}$ is an explicit form of $\epsilon$ in Eq. (\ref{eq:suboptimal_setvalue_mapping}) and $\hat{e}_{y}$ is an upper bound on the suboptimal perturbation. These perturbations can be formulated as a continuous perturbation $\rho\mathbb{B}$, where $\rho = \max\{e_{\tau},e_{y} \}.$ For the upcoming analysis, the next assumption guarantees perturbed solutions to the hybrid system are not arbitrarily far away.
\begin{assumption}\label{assum:compactness}
    The set $Y$ is  compact.
\end{assumption}
The next result formalizes stability of solutions to $\mathcal{H}_{\text{MIQP}}$ with sufficiently small perturbations due to suboptimal solver parameters.

\begin{theorem}\label{thm:perturb}
Let Assumption~\ref{assum:warm_start},~\ref{assum:stable_control} and \ref{assum:compactness} hold. Suppose warm starting $\mathcal{H}_{\text{MIQP}}$ with an initial feasible guess and sufficiently small perturbations from suboptimal solver parameters, then any perturbed solution to the set $\mathcal{A}$ can be bounded as
\[
\vert \pi(t,j) \vert _{\mathcal{A}} \leq \beta(\vert \pi(0,0) \vert_{\mathcal{A}} ,t+j)+\sigma\quad \forall (t,j)\in dom\pi,
\]
where $\sigma>0$ is the perturbation on solutions to $\mathcal{H}_{\text{MIQP}}$ relative to the set $\mathcal{A}$ due to the suboptimal solver parameters.

\end{theorem}
\begin{IEEEproof} The same conditions from Lemma~\ref{lem:asym_stable} hold, so for $\phi_{\epsilon}-V=0$ there is uniform global pre-asymptotic stability. When $\phi_{\epsilon}-V\neq0,$ the rest of this proof obtains a robust version for sufficiently small perturbations. Let Assumption~\ref{assum:compactness} hold, then the hybrid system $\mathcal{H}_{e}$ is the same as the hybrid system with disturbances in \cite{sanfelice2021hybrid}[Section 3.3.3].  This means it can be formulated as a hybrid system with a continuous perturbation function $\rho.$ For every flow time $t\geq 0,$ the suboptimal errors can be encoded by $e_{y}(t),e_{\tau}(t) \in \rho \mathbb{B}$. The function $\vert \pi(t,j) \vert_{\mathcal{A}}$ is the Euclidean norm and is a proper indicator for the stabilizing set $\mathcal{A}$. Then all conditions from Proposition~\ref{prop:robust_KL} hold and the hybrid system $\mathcal{H}_{\text{MIQP}}$ is semiglobally practically robust asymptotic stable, which is equivalent to a $\mathcal{KL}$ function defined as
\[
\vert \pi(t,j) \vert _{\mathcal{A}} \leq \beta(\vert \pi(0,0) \vert_{\mathcal{A}} ,t+j)+\sigma\quad \forall (t,j)\in dom\pi,
\]
with the perturbation term $\sigma>0.$ This completes the proof.
\end{IEEEproof} 

The above theorem provides a hybrid systems based error bound on suboptimal MIQP feedback control solutions due to suboptimal solver parameters.
Theorem~\ref{thm:perturb} is for constant suboptimal solver parameter perturbations. If a practitioner wants a time-varying perturbation, then the next section is relevant.

\subsection{Impact of Time-Varying Suboptimal Solver Parameters}
In Theorem~\ref{thm:perturb}, the suboptimal solver parameters act as a constant perturbation $\sigma$ on solutions to the hybrid system $\mathcal{H}_{\text{MIQP}}$ and is akin to a classic unknown disturbance term. However, unlike classic disturbances, a practitioner can tune a time varying suboptimal solver parameter to balance compute and suboptimality. For illustration purposes, take the case of increasing the maximum number of solver iterates every sample time, then when the solver parameters are the only perturbation on $\mathcal{H}_{\text{MIQP}}$ the term $\sigma$ goes to zero for some finite time $P$, i.e. $\lim_{t\to P}\sigma(t)=0.$ It is beyond the scope of this work to dive further into this notion, but it is a promising theoretical avenue to leverage solver parameters in the same way as step size selection for gradient descent.

Theorem~\ref{thm:perturb} has strong parallels to suboptimal continuous optimization for model predictive control results in the literature \cite{behrendt2023autonomous,liao2020time}. Due to these strong parallels our simulations focus on MIQP variants of model predictive control, which is a common formulation for mixed-logical dynamical systems and signal temporal logic \cite{bemporad1999control,yang2020continuous}.
\begin{remark}
The hybrid model $\mathcal{H}_{\text{MIQP}}$ can address nonlinear dynamics, but the simulation focuses on linear dynamics because Problem~\ref{prob:MIQP} has affine constraints. For a system with nonlinear dynamics $f(x,\kappa(y)),$ a practitioner can linearize $f(x,\kappa(y))$ about a reference trajectory \cite{berberich2022linear}.
\end{remark}

\section{Simulations}\label{sec:simulations}This section analyzes three suboptimal solver parameters: soft time limit, soft memory limit, and iteration limit for MIQP model predictive control (MPC) on the rendezvous and proximity spacecraft problem. These simulations tie into a growing body of literature analyzing abstract models for mixed-integer solver parts \cite{kazachkov2024abstract,le2017abstract} in simulation, where this work extends the abstract modeling to mixed-integer feedback control.

Analyzing soft time limit, soft memory limit, and iteration limit highlight existing solver parameters that are viable under Theorem~\ref{thm:perturb} and useful for practitioners. These solver parameters are explained briefly in Section~\ref{sec:solver_params} and in further detail in \cite{gurobi_MIP}. An important distinction is that soft time limit and soft memory limit produce inexact solutions, while iteration limit with warm starting produces a feasible solution.

The solver parameters are time invariant for each simulation run. All simulations are in MATLAB and YALMIP with Gurobi as the MIQP solver. For the MPC details, a horizon of $N=15$ is used with a sampling time, $\Delta t$, of 5 minutes for slow spacecraft dynamics, see Appendix~\ref{appen:rend_dyn} for more details on the dynamics.

\subsection{Time-Varying MPC}
In Example~\ref{prob:limit_MIQP}, the constraint $\eta(k) = x_{k}$ and objective vector$~\eta(x(\tau))=c$ are time varying and the rest of the constraints are time invariant, the decision vector is defined as $y=[\zeta^{T}~v^{T}]^{T}$, which takes the form of Problem~\ref{prob:MIQP}. Example~\ref{prob:limit_MIQP} has an $\ell^1$ norm control constraint that is nonconvex and is reformulated into a mixed-integer form, which is explained further in Appendix~\ref{appen:L1norm}. The  $\ell^1$ norm control constrained MPC has the following form 
\begin{example}\label{prob:limit_MIQP}
    \begin{align*}
    &\underset{y}{\text{minimize}}  \enskip y^{T}Qy+c^{T}y\\
    &\text{subject to} 
    \enskip \zeta = x_{k}\\  
     &\qquad \zeta_{k+1} = A\zeta_{k}+Bv_{k}, \enskip \forall k = 1,...,N-1\\
     &\qquad  v_{\min}\leq \vert \vert v_{k}\vert \vert_{1} \leq v_{\max} \enskip \forall k = 1,...,N-1\\
     &\qquad \zeta_{N}\in Z_{\text{terminal}}
    \end{align*}
where $\zeta_{k} \in \mathbb{R}^{q}$ is the dynamic state with sample time $k$,$~v_{k}\in \mathbb{R}^{r}$ is the control at sample time $k$,$~\zeta:=[\zeta^{T}_{1},\zeta^{T}_{2},...\zeta^{T}_{N}]^{T},~v:=[v^{T}_{1},v^{T}_{2},...v^{T}_{N-1}]^{T},~y:=[\zeta^{T}~v^{T}]^{T}\in \mathbb{R}^{n},~A\in \mathbb{Q}^{q\times q}$ are the rendezvous dynamics defined in further detail in Appendix~\ref{appen:rend_dyn}, $B\in \mathbb{R}^{q\times r}$ is an electric control thrust defined in further detail in Appendix~\ref{appen:rend_dyn}, $c\in \mathbb{R}^{n}$ is a time-varying vector constructed with normally distributed random numbers, $Q\in \mathbb{R}^{n\times n}$ is a time-invariant cost on state and control, $v_{\min}\in \mathbb{R},~v_{\max}\in \mathbb{R}$ are lower and upper bounds on control $v$,  and $Z_{\text{terminal}}$ is a compact terminal set.
\end{example}

\begin{figure}[!ht]
\centering

    \includegraphics[scale=0.40]{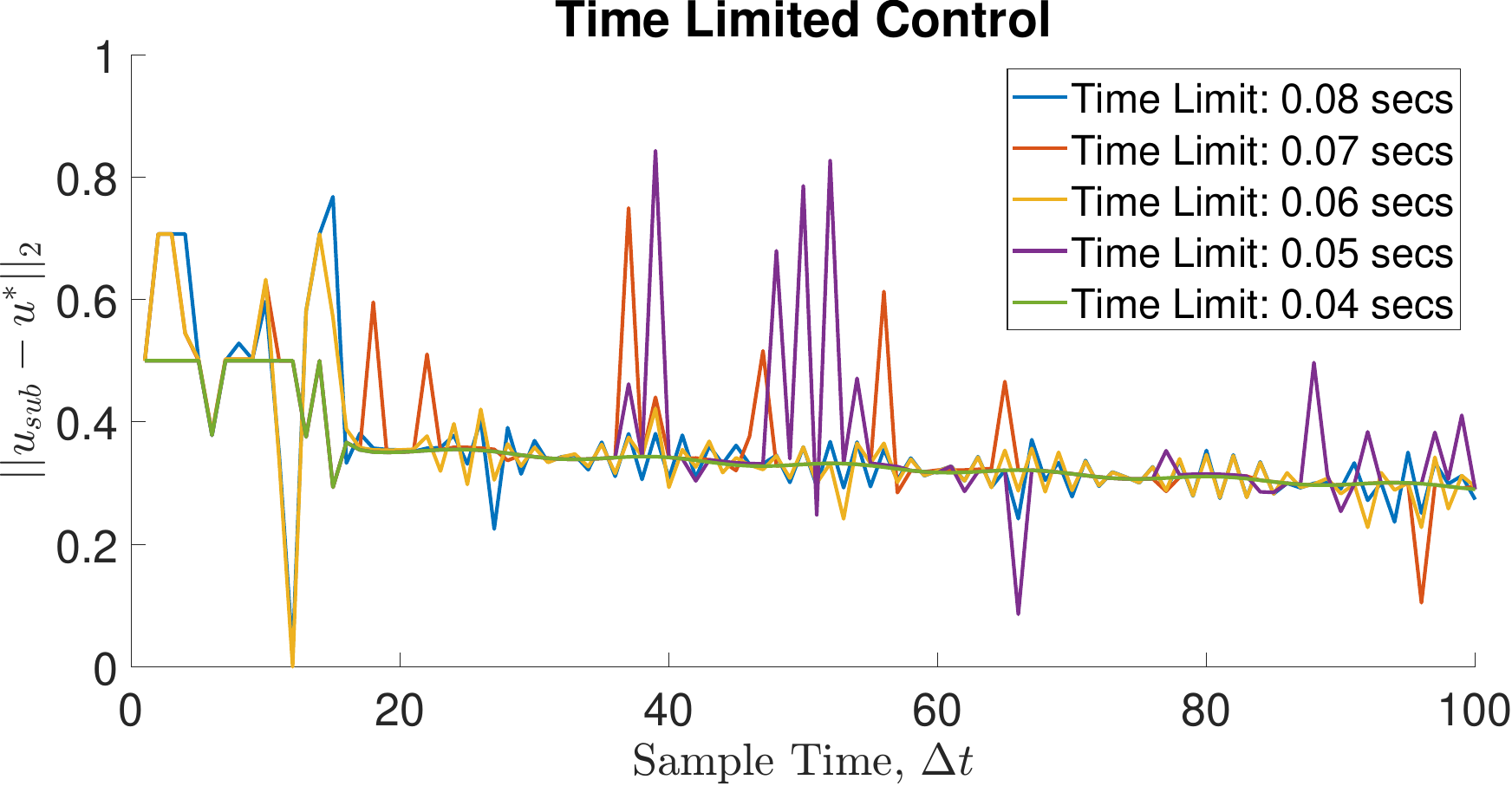}

\caption{Five runs of soft time limited control for Example~\ref{prob:limit_MIQP} are compared. Mild suboptimality was observed for all soft time limits.}
\label{fig:time_control_MPC}
\end{figure}

\begin{figure}[!ht]
\centering

    \includegraphics[scale=0.40]{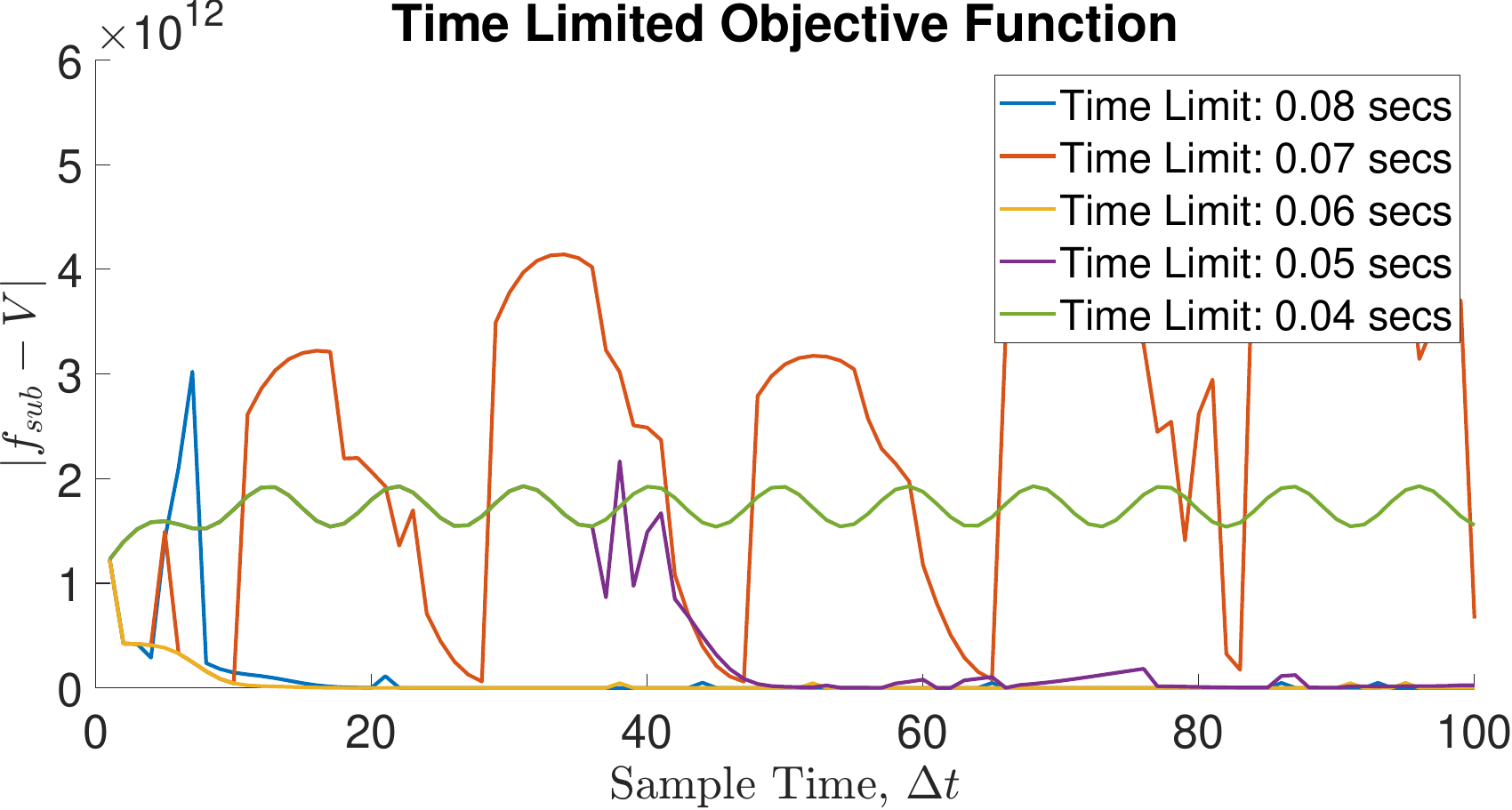}

\caption{Five runs of soft time limited cost function for Example~\ref{prob:limit_MIQP} are compared. Instability is observed when $|f_{sub}-V|$ diverges.}
\label{fig:time_cost_MPC}
\end{figure}

It was observed in Figure~\ref{fig:time_control_MPC} that the suboptimality is mild for soft time limited control. As expected in Figure~\ref{fig:time_cost_MPC}, the soft time limited objective function is stable for soft time limits of $0.08,~0.06,0.05$ seconds and unstable for a soft time limit of 0.04 seconds, which is due to the solver not making sufficient progress and has also been observed in the literature for the quadratic programming MPC case \cite{behrendt2023autonomous}. Unexpectedly in Figure~\ref{fig:time_cost_MPC}, a soft time limit of 0.07 seconds is unstable, which is suspected to be due to poor numerics and deserves further investigation in future works. 

\begin{figure}[!ht]
\centering

    \includegraphics[scale=0.40]{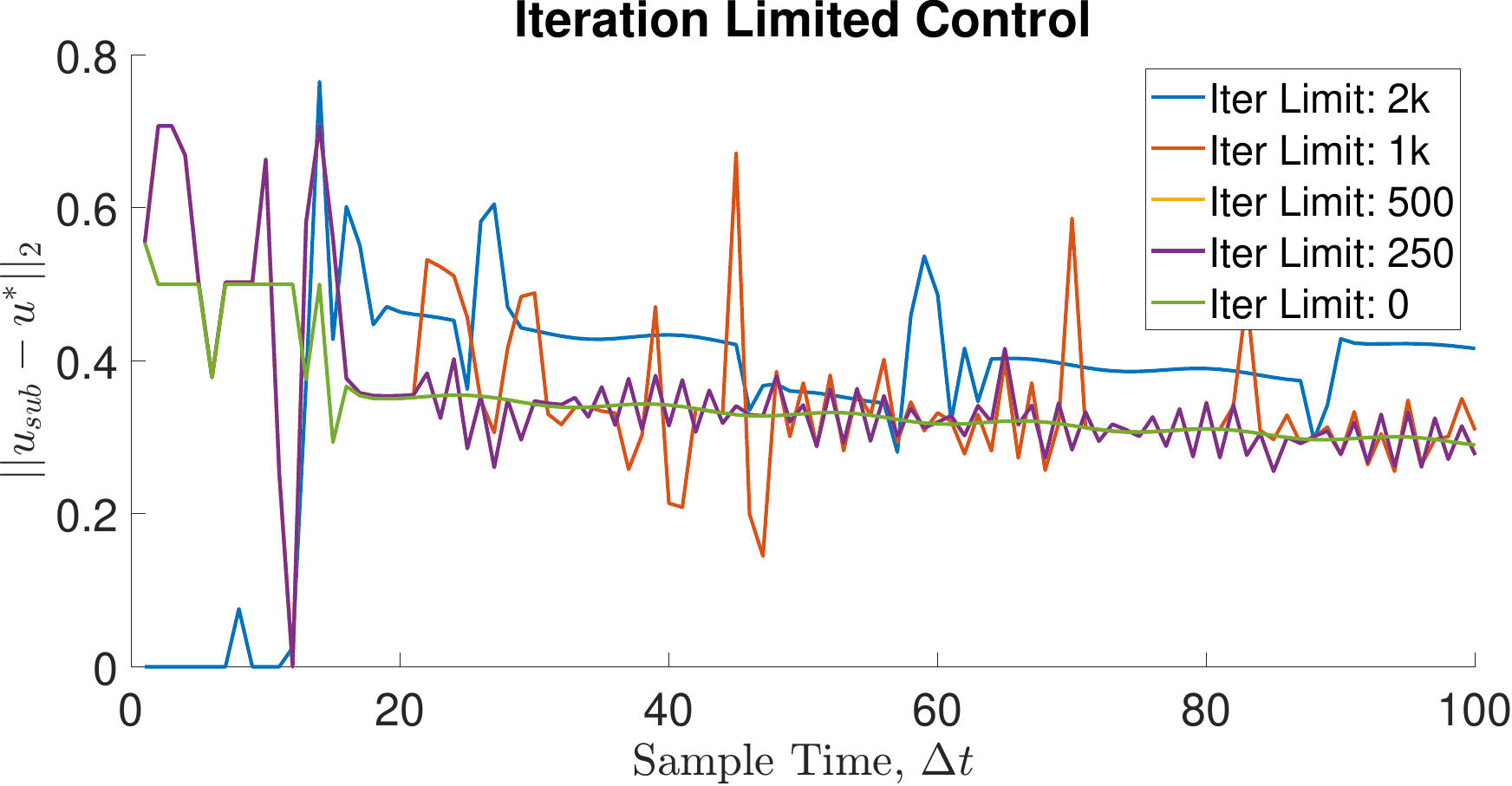}

\caption{Five runs of iteration limited control for Example~\ref{prob:limit_MIQP} are compared. Mild suboptimality was observed for all iteration limits.}
\label{fig:iter_control_MPC}
\end{figure}

\begin{figure}[!ht]
\centering

    \includegraphics[scale=0.40]{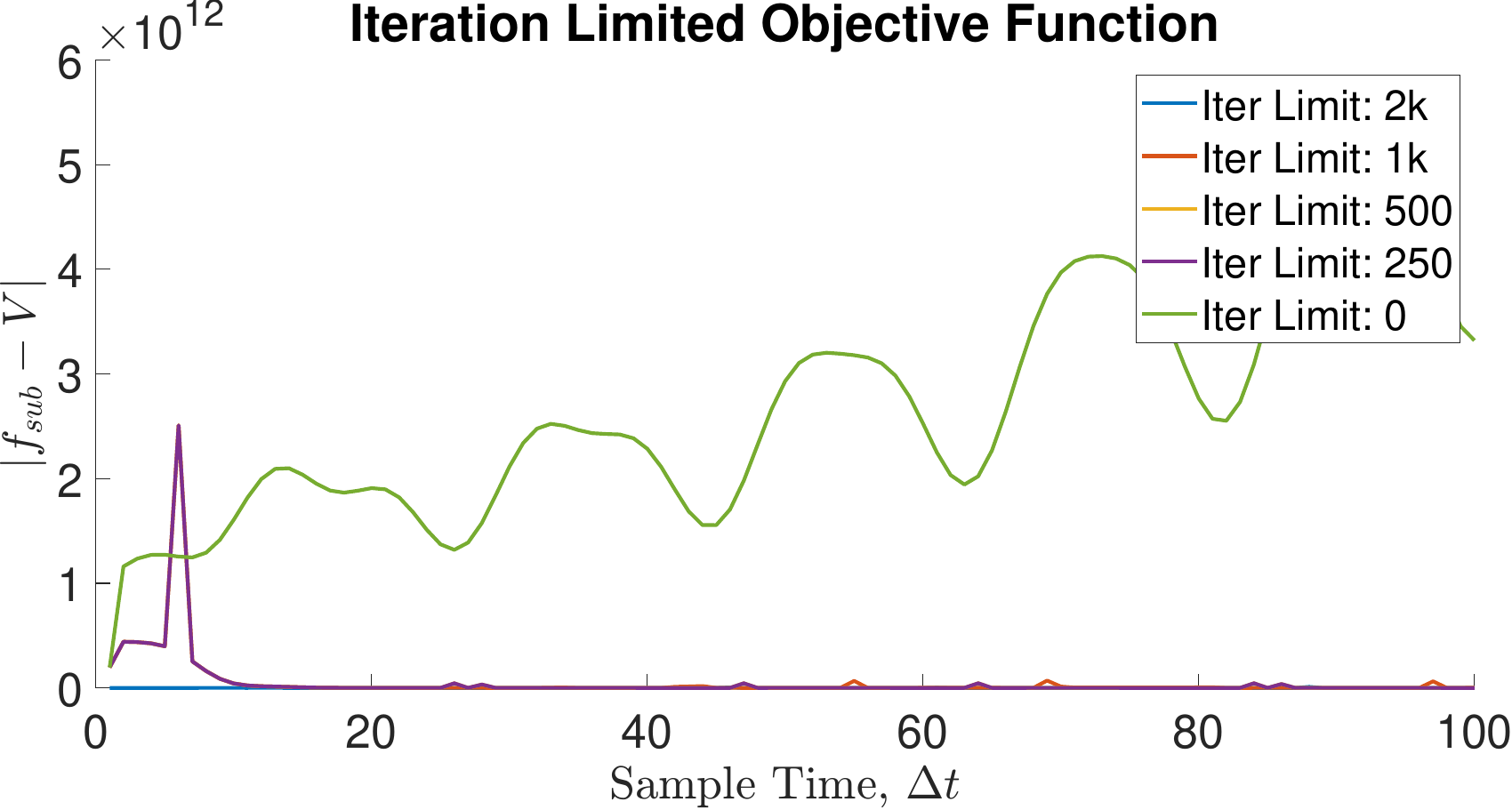}

\caption{Five runs of iteration limited cost function for Example~\ref{prob:limit_MIQP} are compared. Instability is observed when $|f_{sub}-V|$ diverges.}
\label{fig:iter_cost_MPC}
\end{figure}

It was observed in Figure~\ref{fig:iter_control_MPC} that the suboptimality is mild for iteration limited control. As expected in Figure~\ref{fig:iter_cost_MPC}, the iteration limited objective function is stable for all nonzero iteration limits, which is due to a combination of solver heuristics and warm starting.

\begin{figure}[!ht]
\centering

    \includegraphics[scale=0.40]{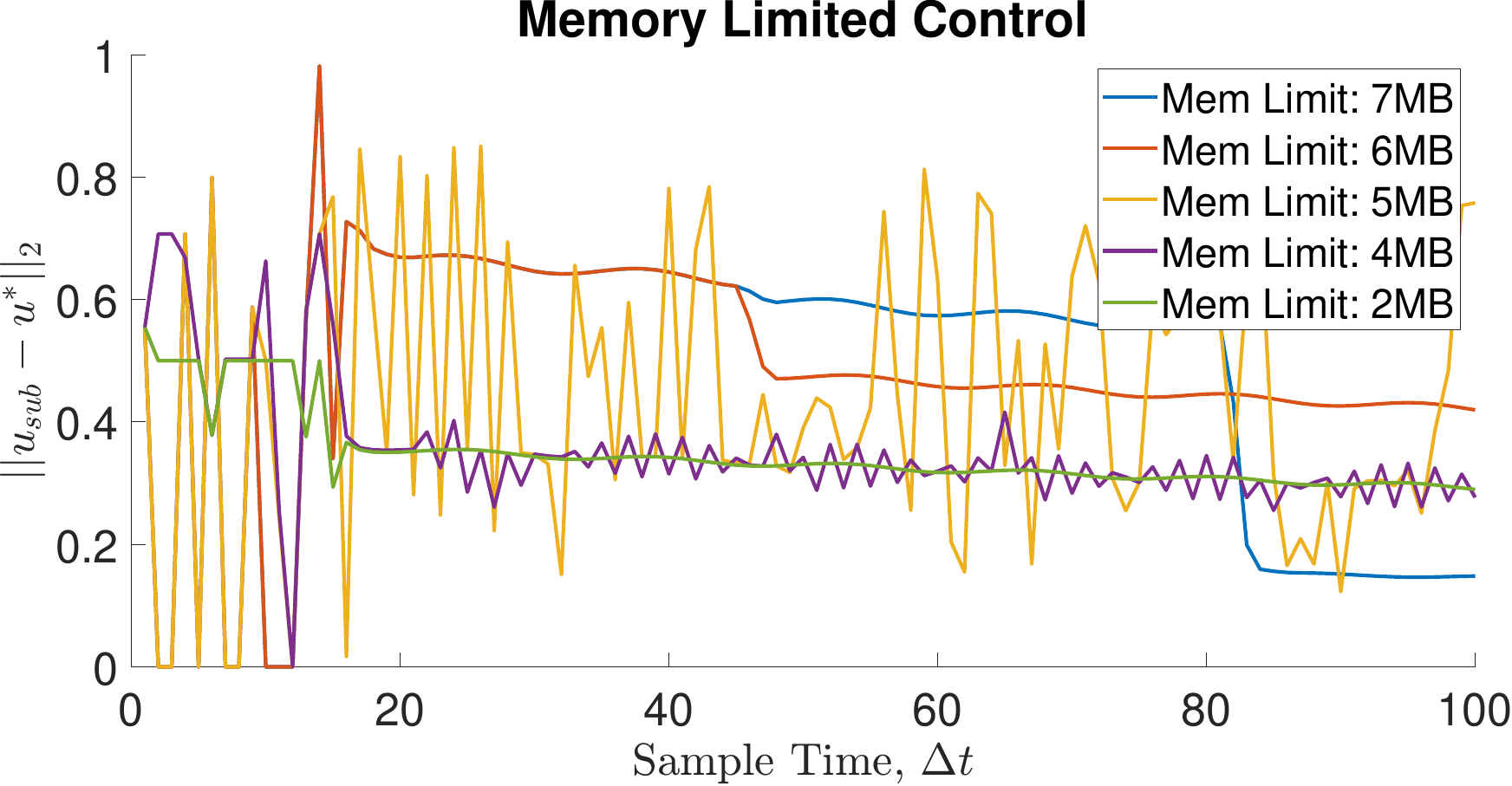}

\caption{Five runs of a soft memory limited control for Example~\ref{prob:limit_MIQP} are compared. Mild suboptimality was observed for all soft memory limits.}
\label{fig:mem_control_MPC}
\end{figure}

It was observed in Figure~\ref{fig:mem_control_MPC} that the suboptimality is mild for soft memory limited control. As expected in Figure~\ref{fig:time_cost_MPC}, the soft memory limited objective function is stable for memory limits of $7$MB,~$6$MB,~$5$MB,~and $4$MB. Additionally, the soft memory limited system is unstable for a soft memory limit of $2$MB, which is likely due to insufficient solver progress.

 Heavy compute in practice is a known concern for MIQPs and in the literature heuristics have been used to obtain a fast solution \cite{takapoui2020simple}, especially for embedded systems. As observed, soft time limit, soft memory limit, and iteration limit parameters are a promising alternative to reduce solver compute for MIQP feedback control while leveraging efficient existing solvers.
 


\begin{figure}[!ht]
\centering

    \includegraphics[scale=0.40]{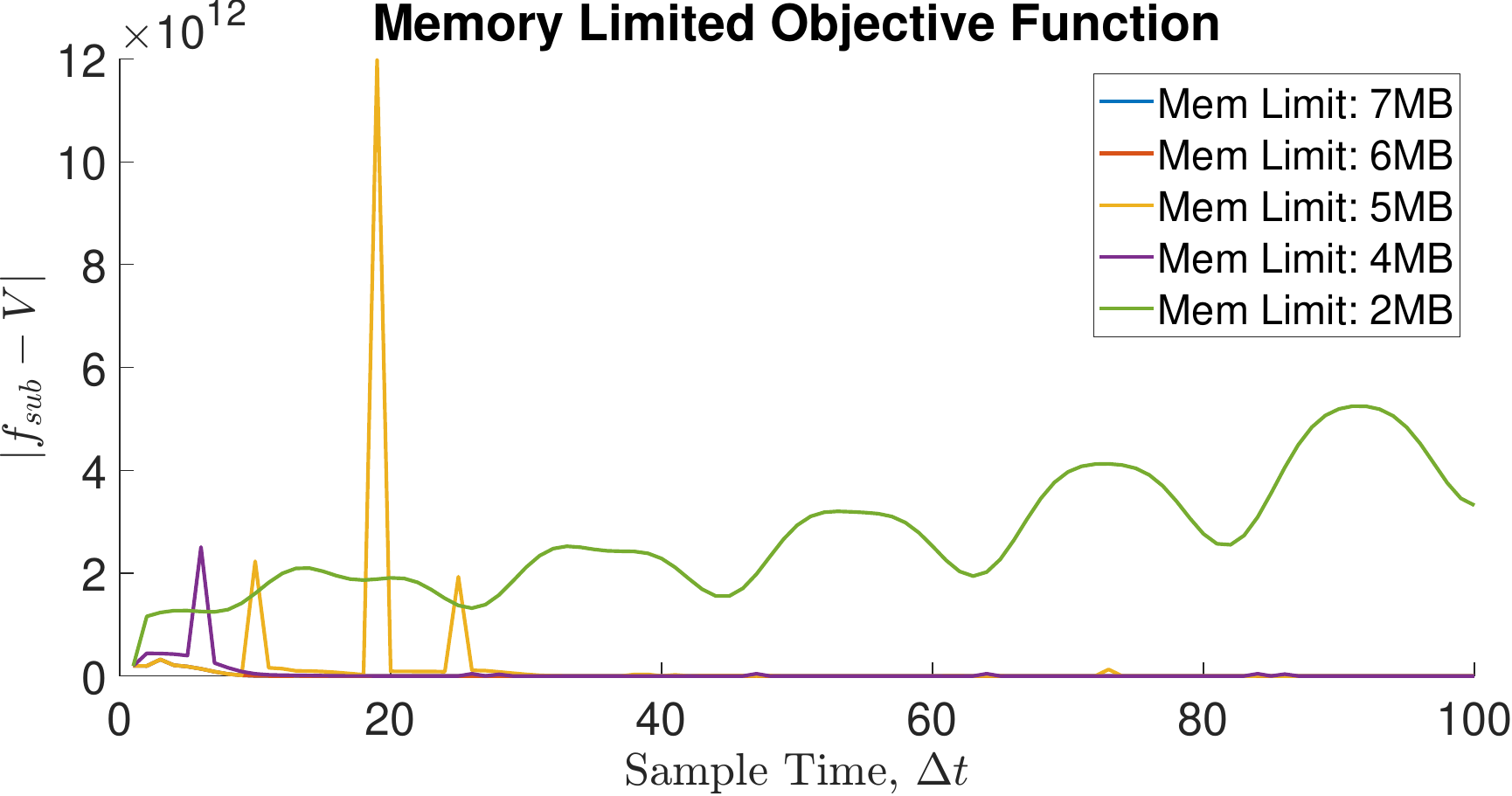}

\caption{Five runs of a soft memory limited cost function for Example~\ref{prob:limit_MIQP} are compared. Instability is observed when $|f_{sub}-V|$ diverges.}
\label{fig:mem_cost_MPC}
\end{figure}

\section{Conclusion}\label{sec:conclusion}
 This paper proposes a hybrid system framework for recursive mixed-integer quadratic programs (MIQPs) controlling a dynamical system. Parametric MIQPs and hybrid systems theory are leveraged to analyze the recursive MIQP and dynamical system together. The MIQP is abstracted to a black box solver with tunable solver parameters for analysis and simulations explore the impact of suboptimal solver parameters on the proposed hybrid system's solutions. Future works will analyze time varying solver parameters and analyze the existence of hybrid solutions with existence solver parameters.


\bibliographystyle{unsrt}
\bibliography{ref}

\section{Appendix}\label{sec:appendix}
This appendix provides technical lemmas, propositions, and definitions for the hybrid system, MIQP stability, and simulation results.
\subsection{Hybrid System Appendix}
This section provides hybrid system related lemmas, propositions, and definitions for completeness. Theorem~\ref{thm:hybrid_cond_MIQP} refers to the following definition, because the definition provides technical details defining upper semi-continuity for set-valued maps. 
\begin{definition}\label{def:u_s_c} \cite{sanfelice2021hybrid} The set-valued map $F:\mathbb{R}^{p}\rightrightarrows\mathbb{R}^{h}$ is upper semi-continuous at $x$ if for every $\nu>0$ there exists $\delta>0$ such that $\zeta \in x +\delta \mathbb{B}$ implies $F(\zeta)\subset F(x)+\nu\mathbb{B}.$ The map $F$ is said to be upper semi-continuous if it upper semi-continuous at each $x\in \mathbb{R}^{p}.$
\end{definition}

Theorem~\ref{thm:perturb} refers to the following Proposition, because it provides technical conditions for stability of solutions to a hybrid system with small perturbations.

\begin{proposition}\label{prop:robust_KL}\cite{sanfelice2021hybrid}[Theorem 3.26] Give a hybrid system $\mathcal{H}$ and a nonempty set $\mathcal{A}\subset \mathbb{R}^{h},$ suppose
\begin{enumerate}
    \item $\mathcal{H}$ satisfies the basic hybrid conditions defined in Definition~\ref{def:hybrid_conds}
    \item $\mathcal{A}$ is a compact pre-asymptotically stable set for $\mathcal{H}$ with basin pre-attraction $\mathcal{B}_{\mathcal{A}}$ (defined in Definition~\ref{def:pre_attractive_basin}), for which the proper indicator $\omega$ of $\mathcal{A}$ on $\mathcal{B}_{\mathcal{A}}$ and the $\mathcal{KL}$ function $\beta$ are such that the $\mathcal{KL}$ function derived from Lemma~\ref{lem:asym_stable} holds for each solution of $\pi$ to $\mathcal{H};$
    \item The $\rho$-perturbation of $\mathcal{H}$ given by $\mathcal{H}_{\rho}$ has a maximum allowable perturbation that is defined example by \cite{sanfelice2021hybrid}[Assumption 3.25].
\end{enumerate} 
Then, the compact set $\mathcal{A}$ satisfies the following $\mathcal{KL}$ function 
\[
\omega(\pi_\rho(t,j))\leq \beta(\omega(y_{\rho}(0,0),t+j)+\kappa \quad \forall(t,j)~dom \pi_{\rho}.
\]
\end{proposition}

Proposition~\ref{prop:robust_KL} refers to the following definition which provides technical details on a basin of pre-attraction for hybrid system.
\begin{definition}\cite{goebel2009hybrid}[Definition 7.3]\label{def:pre_attractive_basin}
Let $\mathcal{H}$ be a hybrid system on $\mathbb{R}^{h}$ and $\mathcal{A}\subset \mathbb{R}^{h}$ be locally pre-asymptotically stable for $\mathcal{H}.$ The basin of pre-attraction of $\mathcal{A}, $ denoted $\mathcal{B}_{\mathcal{A}},$ is the set of points $\zeta\in \mathbb{R}^{h}$ such that every solution $\pi$ to $\mathcal{H}$ with $\pi(0,0)=\zeta$ is bounded and, if it is complete, then also $\lim_{t+j\to \infty}  \vert \pi(t,j) \vert_{\mathcal{A}}=0.$
\end{definition}

\subsection{MIQP Stability Appendix}\label{appen:MIQP_stability}
This appendix builds up technical conditions under Assumption~\ref{assum:rational} where the set value mapping for the value function, constraint set mapping, and $\epsilon$-optimal set value mapping are either continuous or upper semi-continuous, which is crucial to solidify the connection to the basic hybrid conditions in Theorem~\ref{thm:hybrid_cond_MIQP}.
This section starts with a condition for a solution to exist for the MIQP under Assumption~\ref{assum:rational}.
\begin{proposition}\label{prop:thm2.2}
    If Assumption~\ref{assum:rational} holds, then a feasible point exists at which the infimum of the mixed-integer quadratic program is attained if the infimum is finite.
\end{proposition}
\begin{IEEEproof}
See \cite{bank1984stability}[Theorem 2.2].
\end{IEEEproof}
\subsubsection{Continuity of Constraint Mapping}

To discuss the continuity of the constraint set mapping $H: \mathbb{R}^{m}\rightrightarrows\mathbb{R}.$ Let's define the feasible region $H(\eta(k))$ of a perturbed mixed-integer quadratic program on the feasible parameter set as
\[
\Gamma=\{\eta(k)\in \mathbb{R}^{m}\vert H(\eta(k))\neq \emptyset\}.
\]

Upper semi-continuity of Eq. (\ref{eq:suboptimal_setvalue_mapping}) relies on continuity of the constraint mappings $H.$ The following lemma and propositions provide technical details for continuity of the constraint set mapping $H.$

\begin{lemma}\label{lem:cont_contraintmapping}Let the constraint matrix be rational, then the constraint set value mapping is continuous on the feasible parameter set $H(\eta(k)).$
\end{lemma}
\begin{IEEEproof}
    Under Assumption~\ref{assum:rational} the constraint matrix A is rational, then the constraint set mapping $H$ restricted to the feasible parameter set is upper semi-continuous (Proposition~\ref{prop:usc_H}) and the constraint set mapping $H$ is lower semi-continuous (Proposition~\ref{prop:lsc_H}). By definition, a set value mapping is continuous if it is both upper semi-continuous and lower semi-continuous. Both hold, therefore the constraint set mapping is continuous over the feasible parameters. This completes the proof.
\end{IEEEproof}
The following propositions provide technical conditions for lower semi-continuity and upper semi-continuous to hold for the constraint set mapping $H$. 

\begin{proposition}\cite{bank1984stability}\label{prop:usc_H}
If the constraint matrix is rational, then
\begin{enumerate}
    \item there exists a upper semi-continuous set value mapping $K:\Gamma\rightrightarrows\mathbb{R}^{n}$ with compact nonempty images such that the feasible region $H(\eta(k))$ may be represented by 
    \[
    G(\eta(k))=K(\eta(k))+G(0), \eta(k)\in \Gamma,~\text{and}
    \]
    \item the restriction of the constraint set mapping H to the feasible parameter set $\eta(k)$ is upper semi-continuous.
\end{enumerate}
\end{proposition}
\begin{proposition}\cite{bank1984stability}\label{prop:lsc_H}
    Let the constraint matrix be rational and let $\eta(k)_{0}$ be an arbitrarily fixed vector belonging to the feasible set $\Gamma,$ then the restriction of the constraint set value mapping $H$ to the set 
    \[
    \hat{\Gamma}(\eta(k)_{0})=\{ \eta(k) \in \Gamma \vert \Pi_{y_{1},..,y_{s}}H(\eta(k)_{0}\},
    \] where $\Pi$ is the projection operator is lower semi-continuous.
\end{proposition}

\subsubsection{Continuity of Value Function}
The upper semi-continuity relies on continuity of the value function $V.$ The following lemma and proposition provide technical details on continuity of the value function.
\begin{lemma}\cite{bank1984stability}\label{lem:continuous_value_fxn}
    Let $E^{*}$ be a subset of $E$ defined in Proposition~\ref{prop:thm1.2} such that $\eta(k)\in \Gamma^{*}(\eta(k)^{0})$ follows from $(c(k),\eta(k)),~(c(k),\eta(k)^{0})\in E^{*},$ then the rationality of matrix A implies that $V:E^{*}\rightrightarrows \mathbb{R}$ is continuous.
\end{lemma}
\begin{IEEEproof}
    See \cite{bank1984stability}.
\end{IEEEproof}

\begin{proposition}\cite{bank1984stability}\label{prop:thm1.2}
    If the constraint matrix has rational elements only, then the restriction of the value function $V$ defined by $V(c(k),\eta(k))= \inf \{ y^{T}Qy+c(k)^{T}y$ $\vert Ay\leq \eta(k), y_{1},...,y_{s} ~\text{integer}\}$ to the set 
    \[
    E=\{(c(k),\eta(k))\in \mathbb{R}^{hn}\times \mathbb{R}^{m} \vert -\infty <V(c(k),\eta(k))<\infty\}
    \]
    is lower semi-continuous.
\end{proposition}

\subsubsection{Upper Semi-continuity of Suboptimal Mapping}
The following technical proof is provided for completeness as to how upper semi-continuity holds for $\psi_{\epsilon}.$
\begin{IEEEproof}[Proof for Lemma~\ref{lem:subopt_mapping_usc}\cite{bank1984stability}]
    Let $(c(k)^{t},\eta(k)^{t},\epsilon_{t}), $$t= 1,2,...,$ be an arbitrary sequence with 
    \begin{align*}
    (c(k)^{t},\eta(k)^{t})\in E^{*}, t=0,1,2,...,\qquad \lim_{t\to \infty}(c(k)^{t},\eta(k)^{t})=(c(k)^{0},\eta(k)^{0}),\\
    \epsilon_{t}\geq 0, t=1,2,..., \quad \lim_{t\to \infty}\epsilon_{t}=\epsilon_{0.}
    \end{align*}
    Because of the continuity of $V$ on $E^{*}$ we can set $a_{t}=V(c(k)^{t},\eta(k)^{t})+\epsilon_{t}\geq 0$ so that Proposition~\ref{prop:thm2.1} immediately proves the theorem.
\end{IEEEproof}

The following proposition is referred to in the above proof and provides additional technical details.
\begin{proposition}\cite{bank1984stability}\label{prop:thm2.1}
Let the constraint matrix be rational and let the sequence $\{c(k)^{t},\eta(k)^{t},a_{t}\}, t= 1,2,...,$ show the properties
\[
(c(k)^{t},\eta(k)^{t},a_{t})\in E_{N}, t= 1,2,.., \quad \lim_{t\to \infty}(c(k)^{t},\eta(k)^{t},a_{t})=(c(k)^{0},\eta(k)^{0},a_{0})\in E_{N}.
\]
Then, for the level sets $N^{t}=N(c(k)^{t},\eta(k)^{t},a_{t})$ the following holds:
\[
\forall \varepsilon>0~ \exists t(\varepsilon) \quad N^{t}\subset U_{\varepsilon}N^{0}, t>t(\varepsilon)^{2}.
\]
\end{proposition}

\subsection{Mixed-Integer Control Constraint Reformulation}\label{appen:L1norm}
The $\ell^{1}$ norm constraint in Example~\ref{prob:limit_MIQP} is reformulated with binary variables \cite{bertsimas1997introduction} where $z_{k}\in \{ 0,1\}^{m}$ as follows
\begin{align*}
v_{k} = v^{+}_{k}-v^{-}_{k}\qquad \text{for all }k=1,...,N-1,\\
0 \leq v^{+}_{k},v^{-}_{k}\qquad \text{for all }k=1,...,N-1,\\
v^{-}_{k}\leq v_{\max}(\textbf{1}-z_{k}) \qquad \text{for all }k=1,...,N-1,\\
v^{+}_{k}\leq v_{\max}z_{k} \qquad \text{for all }k=1,...,N-1,\\
\sum_{i}v^{+}_{i,k}+\sum_{i}v^{-}_{i,k}\leq v_{\max} \qquad \text{for all }k=1,...,N-1,\\
v_{\min}\leq \sum_{i}v^{+}_{i,k}+\sum_{i}v^{-}_{i,k}\qquad \text{for all }k=1,...,N-1,
\end{align*}
where $v_{\max}\in \mathbb{R}^{m},~v_{\min}\in \mathbb{R}^{m},$ \textbf{1} is a vector of ones, and $m$ is the dimension of the control vector.
\subsection{Rendezvous Dynamics}\label{appen:rend_dyn}
The dynamics for the simulations are the Clohessy Wiltshire dynamics for spacecraft rendezvous where $\Delta t$ is the sample time and the dynamics are
\begin{equation}
    A=\STM{t_0}{t_1}=\begin{bmatrix} \Phi_{\rm rr}(\Delta t) & \Phi_{\rm rv}(\Delta t) \\ \Phi_{\rm vr}(\Delta t) & \Phi_{\rm vv}(\Delta t) \end{bmatrix},
\end{equation}

where 

\begin{align*}
       &\Phi_{\rm rr}(t,t_0)= \begin{bmatrix} 4-3w_c & 0 & 0\\  6(w_s-w)&1&0\\0&0&w_c\end{bmatrix}, 
        &\Phi_{\rm rv}(t,t_0)= \frac{1}{\gamma}\begin{bmatrix} w_s & 2(1-w_c) & 0\\ 2(w_c-1)&4w_s-3w&0\\0&0&w_s\end{bmatrix},
        \\&\Phi_{\rm vr}(t,t_0)= \gamma\begin{bmatrix} 3w_s & 0 & 0\\ 6(w_c-1)&0&0\\0&0&-w_s\end{bmatrix}, 
        &\Phi_{\rm vv}(t,t_0)= \begin{bmatrix} w_c & 2w_s & 0\\ -2w_s&4w_c-3&0\\0&0&w_c\end{bmatrix},        
\end{align*}
with $w=\gamma(\Delta t)$, $\cdot_c=\cos(\cdot)$, and $*_s=\sin(\cdot)$ for the sake of brevity and $\gamma$ is the orbital rate. The control matrix in the simulations uses a continuous thrust model for electric thrust which is modeled as
\begin{equation}
B=\int_{0}^{\Delta t} \STM{\theta}{ t_0} \begin{bmatrix}
    \zeros{m} \\ \eye{m}/m_{s}
\end{bmatrix}d\theta,
\end{equation}
 where $m$ is the dimension of the control, $I_{m}$ is a $m\times m$ identity matrix, $0_{m}$ is a matrix of zeros of dimension $m\times m$, and $m_{s}$ is the mass of the spacecraft. For this work, $m_{s}=100$kg and $\gamma=1.13\times10^{-3}s^{-1}$.

\end{document}